# Multivariate medians and measure-symmetrization*

**Richard A. Vitale**[1]

*University of Connecticut*

**Abstract:** We discuss two research areas dealing respectively with (1) a class of multivariate medians and (2) a symmetrization algorithm for probability measures.

## 1. Introduction

Geometric and stochastic ideas interact in a wide variety of ways over both theory and applications. For two connections that do not seem to have been mentioned before, we present here descriptive comments and suggestions for further investigation. The first deals with *multivariate medians*, an active area of research in which Yehuda himself was interested (see [11] with Cun-Hui Zhang). Using the intrinsic volume functionals for convex bodies, we define a class of multivariate medians and show that among them, as special cases, are the $L_1$ median and the Oja–median. The second question deals with a measure-theoretic generalization of the classic Steiner symmetrization technique for convex bodies.

## 2. A class of multivariate medians

### 2.1. Background

In $\mathbb{R}^d$ with its usual algebraic and metric structures, we consider the class $\mathcal{K}^d$ of convex bodies (compact, convex), which is closed under scaling

$$\lambda K = \{\lambda x : x \in K\}$$

and Minkowski addition

$$K + L = \{x + y : x \in K, y \in L\}.$$

A special class of convex bodies are sums of line segments, or *zonotopes*. The $\lambda$-*parallel body* to $K$ is $K + \lambda B$, where $B_d$ is the closed unit ball in $\mathbb{R}^d$. Steiner's formula gives the volume of a typical parallel body:

$$\mathrm{vol}_d(K + \lambda B_d) = \sum_{j=0}^{d} V_{d-j}(K)\,\mathrm{vol}_j(B_j)\,\lambda^j.$$

---

imsart-lnms ver. 2007/04/13 file: lnms5420.tex date: October 31, 2018

The coefficients $V_j(K)$, $j = 0, 1, 2, \ldots, d$ are the so-called *intrinsic volumes* of $K$ ([10]). They arise in a variety of problems and can, for example, be defined quite differently:

$$V_j(K) = \frac{(2\pi)^{j/2} \, E\mathrm{vol}_j\left(Z_{[j,d]}K\right)}{j! \, \mathrm{vol}_j(B)},$$

where $Z_{[j,d]}$, is a $j \times d$ matrix of independent, standard Gaussian random variables ([17]). Some can be identified with familiar geometric functionals:

$$V_0(K) = 1$$
$$V_1(K) = \text{intrinsic width}$$
$$\vdots$$

$$V_{d-1}(K) = 1/2 \cdot \text{surface area of } K$$
$$V_d(K) = d\text{-dimensional volume of } K$$
$$V_j(K) := 0 \quad \text{for} \quad j > d.$$

Here are some specific intrinsic volumes:

$$V_1([0,1]) = 1$$
$$V_1([a_1, b_1] \times [a_2, b_2] \times \cdots \times a_n, b_n]) = \sum_{i=1}^{n}(b_i - a_i)$$

and generally for $j = 1, 2, \ldots, d$:

$$V_j([a_1, b_1] \times [a_2, b_2] \times \cdots \times [a_n, b_n]) = \sum_{i_1 < i_2 < \cdots < i_j} (b_{i_1} - a_{i_1})(b_{i_2} - a_{i_2}) \cdots (b_{i_j} - a_{i_j}).$$

## 2.2. $V_j$-medians

Suppose that a point sample $x_1, x_2, \ldots, x_n \in \mathbb{R}^d$ is given, and a median is sought. We begin by creating, for every $x \in \mathbb{R}^d$, the zonotope

$$Z(x) = \overline{x - x_1} + \overline{x - x_2} + \cdots + \overline{x - x_n}.$$

$Z(x)$ evidently aggregates, as a polytope, the discrepancies between $x$ and the sample points $x_i$. This can be quantified as follows: from our previous discussion, it is evident that the intrinsic volumes are generally measures of size. More specifically, $V_j$ is a volume-like functional of homogeneity degree $j$. Consider the class of associated variational problems: for each $j = 1, 2, \ldots, d$, minimize $V_j(Z(x))$ over $x \in \mathbb{R}^d$. The minimizing point (if it exists) we call the $V_j$–*median* of the sample.

This point of view unifies two well-known medians that, on other grounds, appear to have little in common:

- The $V_1$–median follows from minimizing $\sum_1^n \|x - x_i\|$ and thus coincides with the well-known $L_1$ median (an interesting treatment of computational issues and other properties appears in [11]).

- The $V_d$-median, on the other hand, depends on minimization of $\text{vol}_d(Z(x))$ and thus coincides with the affine–equivariant median of Oja [7].

Further work:

- It would be of interest to investigate the $V_j$ medians for intermediate values of $j$, $2 \le j \le d-1$. For such a $j$, one would seek to minimize

$$\sum_{1 \le i_1 < i_2 < \cdots < i_j \le n} \left|\det\left(x - x_{i_1}, x - x_{i_2}, \ldots, x - x_{i_j}\right)\right|.$$

- The *Wills functional*

$$W(K) = 1 + V_1(K) + V_2(K) + \cdots + V_d(K)$$

is not homogeneous in $K$ but on account of its other remarkable properties (e.g. [15], [16], [17]) would likely be an interesting alternative size measure. The representation

$$W(Z(x)) = \int_{y \in \mathbb{R}^d} e^{-\pi \text{dist}^2(y, Z(x))} dy$$

could be useful.

- One can focus on the inverse, or *polar*, body related to $Z(x)$. One would then seek to *maximize*, for example, its volume, which is proportional to

$$\int_{u \in S^{d-1}} \frac{1}{\sum_1^n |<u, x - x_i>|^d} du.$$

- Finally, it would be of interest to relate these ideas to other instances of zonotopes' use in data analysis, e.g., [6], [13].

## 3. Steiner symmetrization as mass transport

*Steiner symmetrization* is a geometric transformation of convex bodies that has been useful in a variety of problems, notably when an extremizing body is sought for a prescribed functional. *Mass transport*, on the other hand, embraces a number of issues that deal with shifts of mass (i.e., measure) from one location to another (see, for example, [8, 9]). Recently the question has arisen as to whether it is possible to generalize Steiner symmetrization along the lines of mass transport. Here we suggest a formulation and sketch some preliminary thoughts.

Classical Steiner symmetrization is as follows: given a convex body $K$ and a direction (= unit vector) $u$, locate each chord of $K$ parallel to $u$ and shift it along its line of inclusion so as to re-position its midpoint on the hyperplane $H_u^\perp$ perpendicular to $u$. The aggregate of such shifted chords forms a new convex body, known as the *Steiner symmetral* of $K$.

The original idea for this transformation was apparently due to L'Huiller in the 1780s, but Steiner popularized it in his treatment of the isoperimetric inequality ([3], [4]). Variants of Steiner symmetrization can be found in [1], [3], [5]. The present one is motivated by re-casting the procedure in an equivalent form: regard each chord as bearing a uniformly distributed mass and shift it so as to re-position its center of mass on $H_u^\perp$. More generally one can think of a non-negative (and otherwise nice) $f : \mathbb{R}^d \to \mathbb{R}$, which defines a finite mass distribution on $\mathbb{R}^d$ and therefore on each line $\{u^\perp + tu, -\infty < t < \infty\}$, $u^\perp \in H_u^\perp$. On such a line, the center of mass is

$$\left( \int_{t \in \mathbb{R}} f(u^\perp + tu) dt \right)^{-1} \int_{t \in \mathbb{R}} (u^\perp + tu) f(u^\perp + tu) dt = u^\perp + mu.$$

The shift then amounts to replacing $t \mapsto f(u^\perp + tu)$ with $t \mapsto f(u^\perp + (t-m)u)$. The combined effect of such shifts (over $u^\perp$) amounts to a transport of mass.

One can extend this even farther by regarding total mass as normalized to one and observing that the procedure has an equivalent formulation in terms of conditional expectations of random variables. This avoids issues of regularity for the density. Suppose then that $X$ is a random vector in $\mathbb{R}^d$ for which $EX$ exists ( $\iff E\|X\| < \infty$). Then its *Steiner symmetrization in the direction $u \in S^{d-1}$* is defined to be the random vector

(1) $$X_u = X - E\left[uu'X \mid \Pi_{u^\perp} X\right] = X - uE\left[u'X \mid \Pi_{u^\perp} X\right].$$

Here $\Pi_{u^\perp} = I - uu'$ is (orthogonal) projection onto the subspace $u^\perp$. Thus, conditioned on $\Pi_{u^\perp} X$, the random vector $X$ is shifted so that its conditional expectation lies in $u^\perp$.

It can be shown rather directly that this formulation extends classical Steiner symmetrization:

**Theorem 1.** *Suppose that $X$ is uniformly distributed on a convex body $K$. Then (i) $X_u$ is uniformly distributed on the Steiner symmetral of $K$, and (ii) there is a sequence of symmetrizations that converge in distribution to uniform measure on the (centered) ball of the same volume as $K$.*

As we discuss below, general results appear to be difficult to prove. But the special case of symmetrization of a Gaussian measure already presents some interesting phenomena. We provide details for completeness.

**Theorem 2.** *Suppose that $X$ has Gaussian distribution in $\mathbb{R}^d$ with mean $\mu$ and (invertible) covariance matrix $\Sigma$, $X \sim N(\mu, \Sigma)$. Then (i) any Steiner symmetrization yields another Gaussian distribution, and (ii) there is a sequence of Steiner symmetrizations producing a limiting Gaussian distribution that is spherically symmetric about the origin.*

*Proof.* (i) Direct by the linear nature of the transformation.
(ii) There is no loss of generality in assuming $\mu = 0$ since, if not, a symmetrization with $u = \mu/\|\mu\|$ yields this centering. For the first actual symmetrization, fix $u$ and then as above

$$X_u = X - uE\left[u'X \mid \Pi_{u^\perp} X\right].$$

For convenience, let $\Pi = \Pi_{u^\perp}$. Well-known properties of normal random variables provide that $E[u'X \mid \Pi X]$ is that linear functional of $\Pi X$ which minimizes $E(u'X - w'\Pi X)^2$. Setting $\nabla E[u'X - w'\Pi x]^2 = 0$ leads to

$$\begin{aligned} \Pi \Sigma \Pi' w &= \Pi \Sigma u \\ \iff \quad \Sigma \Pi w &= \Sigma u + cu \quad \text{for a constant } c \\ \iff \quad \Pi w &= u + c\Sigma^{-1} u. \end{aligned}$$

For $c$, we apply $u'$ to get $0 = 1 + cu'\Sigma^{-1}u \Longrightarrow c = -\left(u'\Sigma^{-1}u\right)^{-1}$. Then

$$\Pi w = \Pi^2 w = \Pi\left[u + c\Sigma^{-1}u\right] = c\Pi\Sigma^{-1}u$$

so that

$$E\left[u'X \mid \Pi X\right] = -cu'\Sigma^{-1}\Pi X$$

and

$$\begin{aligned} X_u &= X - cuu'\Sigma^{-1}\Pi X \\ &= (I - cuu'\Sigma^{-1}\Pi)X. \end{aligned}$$

It follows that

$$\begin{aligned} \Sigma_u &= EX_u X_u' \\ &= E(I - cuu'\Sigma^{-1}\Pi)XX'(I - c\Pi\Sigma^{-1}uu') \\ &= (I - cuu'\Sigma^{-1}\Pi)\Sigma(I - c\Pi\Sigma^{-1}uu'). \end{aligned}$$

Now suppose that $(v_1, \lambda_1)$ and $(v_2, \lambda_2)$ are two eigenpairs of $\Sigma$ ($\|v_1\| = \|v_2\| = 1$). A convenient choice for symmetrization is $u = \frac{1}{\sqrt{2}}(v_1 + v_2)$. Then

$$\Pi = I - \frac{1}{2}(v_1 + v_2)(v_1' + v_2'),$$

and

(2) $$\frac{1}{c} = -u'\Sigma^{-1}u = -\frac{1}{2}(v_1 + v_2)\Sigma^{-1}(v_1 + v_2) = -\frac{1}{2}\left(\frac{1}{\lambda_1} + \frac{1}{\lambda_2}\right).$$

Further,

$$\begin{aligned} uu'\Sigma^{-1}\Pi &= \frac{1}{2}(v_1 + v_2)\left(\frac{1}{\lambda_1}v_1' + \frac{1}{\lambda_2}v_2'\right)\left(I - \frac{1}{2}(v_1 + v_2)(v_1' + v_2')\right) \\ &= \frac{1}{2}(v_1 + v_2)\left(\frac{1}{\lambda_1}v_1' + \frac{1}{\lambda_2}v_2'\right) - \frac{1}{4}\left(\frac{1}{\lambda_1} + \frac{1}{\lambda_2}\right)(v_1 + v_2)(v_1' + v_2') \\ &= \frac{1}{4}(v_1 + v_2)\left[2\frac{1}{\lambda_1}v_1' + 2\frac{1}{\lambda_2}v_2' - \frac{1}{\lambda_1}v_1' - \frac{1}{\lambda_1}v_2' - \frac{1}{\lambda_2}v_1' - \frac{1}{\lambda_2}v_2'\right] \\ &= \frac{1}{4}(v_1 + v_2)\left[\frac{1}{\lambda_1}(v_1' - v_2') + \frac{1}{\lambda_2}(v_2' - v_1')\right] \\ &= \frac{1}{4}\left(\frac{1}{\lambda_1} - \frac{1}{\lambda_2}\right)(v_1 + v_2)(v_1' - v_2'). \end{aligned}$$

With (2), this gives

$$cuu'\Sigma^{-1}\Pi = -\frac{1}{2}\frac{\lambda_2 - \lambda_1}{\lambda_1 + \lambda_2}(v_1 + v_2)(v_1' - v_2').$$

Then

$$\Sigma_u = \left(I + \frac{1}{2}\frac{\lambda_2 - \lambda_1}{\lambda_1 + \lambda_2}(v_1 + v_2)(v'_1 - v'_2)\right) \Sigma \left(I + \frac{1}{2}\frac{\lambda_2 - \lambda_1}{\lambda_1 + \lambda_2}(v_1 - v_2)(v'_1 + v'_2)\right)$$

$$= \left(\Sigma + \frac{1}{2}\frac{\lambda_2 - \lambda_1}{\lambda_1 + \lambda_2}(v_1 + v_2)(\lambda_1 v'_1 - \lambda_2 v'_2)\right)$$
$$\times \left(I + \frac{1}{2}\frac{\lambda_2 - \lambda_1}{\lambda_1 + \lambda_2}(v_1 - v_2)(v'_1 + v'_2)\right)$$

$$= \Sigma + \frac{1}{2}\frac{\lambda_2 - \lambda_1}{\lambda_1 + \lambda_2}(v_1 + v_2)(\lambda_1 v'_1 - \lambda_2 v'_2)) + \frac{1}{2}\frac{\lambda_2 - \lambda_1}{\lambda_1 + \lambda_2}(\lambda_1 v_1 - \lambda_2 v_2)(v'_1 + v'_2)$$
$$+ \left[\frac{1}{2}\frac{\lambda_2 - \lambda_1}{\lambda_1 + \lambda_2}\right]^2 (v_1 + v_2)(\lambda_1 v'_1 - \lambda_2 v'_2)(v_1 - v_2)(v'_1 + v'_2)$$

$$= \Sigma + \frac{1}{2}\frac{\lambda_2 - \lambda_1}{\lambda_1 + \lambda_2}\left[(v_1 + v_2)(\lambda_1 v'_1 - \lambda_2 v'_2) + (\lambda_1 v_1 - \lambda_2 v_2)(v'_1 + v'_2)\right.$$
$$\left. + \frac{1}{2}(\lambda_2 - \lambda_1)(v_1 + v_2)(v'_1 + v'_2)\right]$$

$$= \Sigma + \frac{1}{2}\frac{\lambda_2 - \lambda_1}{\lambda_1 + \lambda_2}\left[\frac{3\lambda_1 + \lambda_2}{2}v_1 v'_1 - \frac{\lambda_1 + 3\lambda_2}{2}v_2 v'_2 - \frac{1}{2}(\lambda_2 - \lambda_1)(v_1 v'_2 + v_2 v'_1)\right].$$

Ignoring all but the first two terms of the spectral decomposition $\Sigma = \lambda_1 v_1 v'_1 + \lambda_2 v_2 v'_2 + \cdots + \lambda_d v_d v'_d$, one gets

$$\lambda_1 v_1 v'_1 + \lambda_2 v_2 v'_2 + \frac{1}{2}\frac{\lambda_2 - \lambda_1}{\lambda_1 + \lambda_2}\left[\frac{3\lambda_1 + \lambda_2}{2}v_1 v'_1 - \frac{\lambda_1 + 3\lambda_2}{2}v_2 v'_2 - \frac{1}{2}(\lambda_2 - \lambda_1)(v_1 v'_2 + v_2 v'_1)\right]$$

$$= \frac{\lambda_1^2 + 6\lambda_1\lambda_2 + \lambda_2^2}{4(\lambda_1 + \lambda_2)}(v_1 v'_1 + v_2 v'_2) - \frac{1}{4}\frac{(\lambda_2 - \lambda_1)^2}{\lambda_1 + \lambda_2}(v_1 v'_2 + v_2 v'_1).$$

As a $2 \times 2$ array, this has eigenvalues

$$\frac{\lambda_1^2 + 6\lambda_1\lambda_2 + \lambda_2^2}{4(\lambda_1 + \lambda_2)} \pm \frac{(\lambda_2 - \lambda_1)^2}{4(\lambda_1 + \lambda_2)},$$

which simplify to

$$\frac{1}{2}(\lambda_1 + \lambda_2), \quad \left[\frac{1}{2}\left(\frac{1}{\lambda_1} + \frac{1}{\lambda_2}\right)\right]^{-1}.$$

It follows that $\Sigma_u$ has the same eigenvalues as $\Sigma$ except for $\lambda_1$ and $\lambda_2$, which are replaced by their arithmetic and harmonic means. Iterating this double-mean transformation leads to a common limit ([2]), i.e., $\lambda_1$ and $\lambda_2$ can be replaced by a pair of eigenvalues arbitrarily close to one another. Since this is true for any pair of eigenvalues, the result is established. □

**Remarks**

(1) The appearance of the two means in the proof is of historical interest, as the so–called "double-mean-iteration." See, for example, [2] with a discussion of the work of Archimedes and Gauss among others. It would be of interest to know whether it appears here merely as a curiosity or whether there is a deeper theory connecting symmetrization of Gaussian measures and eigenvalue means.

(2) In analogy with the classical theory, one would like to have results along the following lines:

**Conjecture 1.** *For some sequence of Steiner symmetrizations [resp. for almost every sequence of Steiner symmetrizations (i.e. utilizing independent, random selections of u, cf. [14])], the induced sequence of symmetrized distributions is asymptotically spherically symmetric.*

**Conjecture 2.** Any *sequence of Steiner symmetrizations yields a (weakly) convergent sequence of probability measures.*

These do not seem easy to establish, but it seems plausible that results along the following lines may be helpful:

**Proposition 1.** *Steiner symmetrization is norm-reducing in mean-square:*

$$E\|X\|^2 - E\|X_u\|^2 = E\|X - X_u\|^2$$
$$= E\left[E(u'X|\Pi_{u^\perp}X)\right]^2 \geq 0.$$

*Proof.* The following case is generic: $d = 2$, $u = \binom{1}{0}$

$$X = \begin{pmatrix} X^1 \\ X^2 \end{pmatrix} \longrightarrow X_u = \begin{pmatrix} X^1 - E(X^1|X^2) \\ X^2 \end{pmatrix}$$

$$E\|X\|^2 - E\|X_u\|^2 = E\|X - X_u\|^2$$
$$= E\left[E(X^1|X^2)\right]^2 \geq 0. \qquad \square$$

(3) Finally, it would be of interest to place the transformation $X \longmapsto X_u$ in a martingale context.